\documentclass{amsart}
\usepackage[T1]{fontenc}
\usepackage{textcomp}
\usepackage[utf8x]{inputenc}
\usepackage[english]{babel}
\usepackage{ucs}
\usepackage{graphicx}
\usepackage{mathrsfs}
\usepackage{rotating} 
\usepackage{tikz}
\usepackage{amssymb}
\usepackage[titletoc,toc]{appendix}
\numberwithin{equation}{section} 
\newtheorem{example}{Example}[section]

\theoremstyle{definition}

\theoremstyle{remark}

\numberwithin{equation}{section}

\newcommand{\A}{\mathcal{A}}

\newcommand{\F}{\mathbb{F}}

\newcommand{\C}{\mathbb{C}}

\newcommand{\Sp}[1]{\mathrm{Sp}_2(#1)}

\makeindex
\begin{document}

\title{The equivariant Euler characteristic of~$\mathcal{A}_3[2]$}%
\author{Jonas Bergstr\"om}%
\author{Olof Bergvall}%
\address{Matematiska Institutionen, Stockholms Universitet, 106 91, Stockholm, Sweden}
\email{jonasb@math.su.se}
\address{Matematiska Institutionen, Uppsala Universitet, Box 480, 751 06, Sweden}
\email{olof.bergvall@math.uu.se}

\begin{abstract}
We compute the weighted Euler characteristic, equivariant with respect to the action of the symplectic group of degree six over the field of two elements, of the moduli space of principally polarized abelian threefolds together with a level two structure. 
\end{abstract}

\maketitle

\section{Introduction}
Let $\mathcal{A}_{g}[2]$ denote the moduli space of principally polarized 
abelian varieties of dimension $g$ together with a full level $2$ structure. 
Similarly, let $\mathcal{M}_g[2]$ denote the moduli space of smooth  
curves of genus $g$ together with a full level $2$ structure.
Note that we are considering these as \emph{coarse} moduli spaces. 
The two 
types of moduli spaces are connected through the Torelli 
morphism $t_g: \mathcal{M}_g[2] \to \mathcal{A}_{g}[2]$ sending a 
smooth curve to its Jacobian. There is an action of $\Sp{2g}$,  the symplectic group 
of degree $2g$ over the field of two elements, on $\mathcal{A}_{g}[2]$ and 
$\mathcal{M}_{g}[2]$ via its action on the level $2$ structure (for some more details 
see for instance \cite{bergvallq2}). 

For $g \leq 2$, $\mathcal{M}_{g}[2]$ is isomorphic to the locus $\mathcal{H}_{g}[2]$ consisting of 
hyperelliptic curves. The space $\mathcal{H}_{g}[2]$ can in turn be described as a disjoint union 
of copies of $\mathcal{M}_{0,2g+2}$ (see \cite{dolgachevortland}, \cite{runge} and \cite{tsuyumine}), the moduli space of smooth genus $0$ curves together with 
$2g+2$ marked points. The cohomology of $\mathcal{M}_{0,2g+2}$, 
 together with the action of the symmetric group $S_{2g+2}$, 
can (because of purity) be computed using counts of points over finite fields (see for instance \cite{dimcalehrer} and \cite{kisinlehrer}). 
In Section~\ref{sec-g1} respectively  
Section~\ref{sec-g2} below, we compute in this way the  $\Sp{2}$-action and Hodge structure of the cohomology 
of $\mathcal{M}_{1}[2] \cong \mathcal{A}_{1}[2]$ respectively  the  $\Sp{4}$-action and Hodge structure 
of the cohomology of  $\mathcal{M}_{2}[2]$. By adding the complement of 
$t_2(\mathcal{M}_{2}[2])$ inside $\mathcal{A}_{2}[2]$ consisting of products of elliptic curves we 
also compute the $\Sp{4}$-equivariant weighted Euler characteristic of $\mathcal{A}_{2}[2]$. For a definition of 
this type of Euler characteristic see Section~\ref{sec-euler}. 

The main result of this article is Table~\ref{a3table} which contains the $\Sp{6}$-equivariant weighted Euler characteristic of 
$\mathcal{A}_{3}[2]$. This is based upon the work of the second author in \cite{bergvallq2} in which the 
$\Sp{6}$-equivariant cohomology of $\mathcal{M}_{3}[2]$ is computed, see Section~\ref{sec-decomp} below.  
There are two other loci consisting of either products of an indecomposable abelian surface and an elliptic curve,
or products of three elliptic curves. The cohomology of these loci are computed in Section~\ref{sec-A21} and 
Section~\ref{sec-A111} respectively.

We note in Section~\ref{sec-main} that the weighted Euler characteristic of $\mathcal{A}_{3}[2]$ contains much fewer 
classes than the weighted Euler characteristic of its different loci. This cancellation property was noted also in \cite{jbgvdg} 
for the integer valued Euler characteristic of local systems upon the corresponding strata inside $\mathcal{A}_{3}$,  
the moduli space of principally polarized abelian threefolds with no level structure. The Hodge structure of the 
cohomology of $\mathcal{A}_{3}$ was previously known, see \cite{Hain}.

\section{Euler characteristics} \label{sec-euler}
For a quasi-projective variety $X$ defined over $\C$, 
let $W_{k}H^i(X)$ denote the weight $k$-part of $H^i(X)$, the $i$th Betti cohomology group 
with complex coefficients.  
For an action of a finite group $G$ on $X$, let the 
$G$-equivariant weighted Euler characteristic of $X$ be the virtual representation of $G$ defined as
\begin{equation*}
 e_X(v) = \sum_{i,k \geq 0} (-1)^i \, W_kH^i(X) \, v^k.
\end{equation*}
This Euler characteristic is additive in the sense that if $X=Y \sqcup Z$, where $Y$ and $Z$ are preserved by $G$, 
then $e_{X}(v) =v^{2 \mathrm{codim}_X(Y)} e_Y(v) + v^{2 \mathrm{codim}_X(Z)} e_Z(v)$. Note that if $X$ fullfils 
purity, in the sense of Dimca and Lehrer in~\cite{dimcalehrer}, then one can from this Euler characteristic 
determine the individual cohomology groups as representations of $G$.

Say now that $X$ is a variety defined over $\mathcal{O}[\frac{1}{N}] $, where $\mathcal{O}$ is a ring of integers of 
an algebraic number field, together with an action of a finite group $G$.  Say furthermore that there is a polynomial 
$P(t)$, with complex coefficients and of degree $2\dim X$, such that $P(q)=|X(\F_q)|$ for almost all prime powers $q$.  
We can, using the Lefschetz fixed point theorem, from this information determine the weighted Euler characteristic of $X(\C)$. 
The set $X(\F_q)$ consists of the fixed points of Frobenius. By counting the fixed points of Frobenius 
composed with elements of $G$ we can in the same way determine the $G$-equivariant weighted Euler characteristic 
of $X(\C)$. This will be called a twisted point count. For a reference, see \cite[Appendix~A]{dimcalehrer2}.

\section{Decomposable and indecomposable abelian threefolds} \label{sec-decomp} 
We say that an abelian threefold is indecomposable if it is not isomorphic to
a product of abelian varieties of lower dimension. We denote the corresponding 
locus in $\mathcal{A}_3[2]$ by $\mathcal{A}_3^{\mathrm{in}}[2]$.

The Torelli morphism $t_3$ gives
an isomorphism $\mathcal{M}_3[2] \cong \mathcal{A}_3^{\mathrm{in}}[2]$
(on the level of coarse moduli spaces). The moduli space $\mathcal{M}_3[2]$ can
be decomposed as a disjoint union
\begin{equation*}
 \mathcal{M}_3[2] = \mathcal{Q}[2] \sqcup \mathcal{H}_3[2]
\end{equation*}
where $\mathcal{Q}[2]$ denotes the locus consisting of curves whose canonical model
is a plane quartic curve and where $\mathcal{H}_3[2]$ 
denotes the hyperelliptic locus.

The cohomology groups of $\mathcal{Q}[2]$ and $\mathcal{H}_3[2]$ were
determined as representations of $\Sp{6}$ by the second author in \cite{bergvallq2}.
For completeness, we repeat the results in Table~\ref{Qtable} and Table~\ref{Hyptable}. 

There are two types of decomposable abelian threefolds. 
The threefold can  either be isomorphic to a product of an indecomposable abelian surface and an elliptic curve 
or to a product of three elliptic curves.
We denote the corresponding loci in $\mathcal{A}_3[2]$ by $\mathcal{A}_{2,1}[2]$
and $\mathcal{A}_{1,1,1}[2]$ respectively.

\section{The main result} \label{sec-main} 
We have the decomposition 
\begin{equation*}
 \mathcal{A}_3[2] = 
 t_3(\mathcal{Q}[2]) \sqcup
 t_3(\mathcal{H}_3[2]) \sqcup
 \mathcal{A}_{2,1}[2] \sqcup
 \mathcal{A}_{1,1,1}[2]
\end{equation*}
and below we compute the cohomology groups of each of the spaces on the right hand side as representations of $\Sp{6}$. 
Moreover, we will see that each cohomology group $H^i$ of a space on the right hand side is pure of weight $2i$ 
and Tate type $(i,i)$. 

By the additivity of the weighted Euler characteristic, 
\begin{equation*}
 e_{\mathcal{A}_3[2]}(v) = 
 e_{t_3(\mathcal{Q}[2])}(v) +
 v^2 \, e_{t_3(\mathcal{H}_3[2])}(v) +
 v^4 \, e_{\mathcal{A}_{2,1}[2]}(v) +
 v^6 \, e_{\mathcal{A}_{1,1,1}}(v).
\end{equation*}
Putting the results together for the different strata we get the $\Sp{6}$-equivariant weighted Euler characteristic of 
$\mathcal{A}_3[2]$, see Table~\ref{a3table}. Each column in this table corresponds to an irreducible representation of $\Sp{6}$.
The irreducible representations are denoted $\phi_{dn}$ where $d$ is the dimension of
the representation and $n$ is letter used to distinguish different representations of the same
dimension, see \cite{conwayetal}.

\begin{table}[ht]
\begin{equation*}
\resizebox{0.9\textwidth}{!}{$
 \begin{array}{c|cccccccccc}
  \, & \phi_{1a} & \phi_{7a} & \phi_{15a} & \phi_{21a} & \phi_{21b} & \phi_{27a} & \phi_{35a} & \phi_{35b} & \phi_{56a} & \phi_{70a} \\
  \hline
   e_{\mathcal{A}_3[2]}(v) & 1+v^2+v^4+v^6+v^{12} & 0 & v^{12} & 0 & 0 & -v^6-v^8 & 0 & -v^6-v^8+v^{12} & 0 & 0 \\
   \hline
 \, & \phi_{84a} & \phi_{105a} & \phi_{105b} & \phi_{105c} & \phi_{120a} & \phi_{168a} & \phi_{189a} & \phi_{189b} & \phi_{189c} & \phi_{210a} \\
 \hline
   e_{\mathcal{A}_3[2]}(v) & v^{12} & 0 & v^4 & 0 & v^{10} & v^{10} & 0 & 0 & 0 & -v^6 \\
   \hline
 \, & \phi_{210b} & \phi_{216a} & \phi_{280a} & \phi_{280b} & \phi_{315a} & \phi_{336a} & \phi_{378a} & \phi_{405a} & \phi_{420a} & \phi_{512a} \\
 \hline
   e_{\mathcal{A}_3[2]}(v) & -v^6 &  0 &  0 &  v^{10} &  0 &  0 & 0 & 0 & v^8 & -v^{12}
 \end{array}$}
\end{equation*}
\caption{The $\Sp{6}$-equivariant weighted Euler characteristic of $\mathcal{A}_3[2]$.}
\label{a3table}
\end{table}

Note that the results for $\phi_{1a}$ agree with the computation of the cohomology groups of $\mathcal{A}_3$ together with their Hodge structure in \cite{Hain}. Note also that only 13 of the 30 irreducible representations of $\Sp{6}$ occur in $e_{\mathcal{A}_3[2]}(v)$
and that for each irreducible representation the coefficients of $v^i$ are all either zero or $\pm 1$.
This is in sharp contrast to the cohomology of the individual pieces - all irreducible representations
except $\phi_{7a}$ occur in some cohomology group of some piece and they occur with multiplicities
up to $14$.

\section{An indecomposable abelian surface and an elliptic curve} \label{sec-A21} 
As in the genus $3$ case, $t_2$ gives an isomorphism $\mathcal{M}_2[2] \cong \mathcal{A}_2^{\mathrm{in}}[2]$, 
where $ \mathcal{A}_2^{\mathrm{in}}[2]$ denotes the indecomposable locus inside $ \mathcal{A}_2[2]$. 

There is a close relationship between $\mathcal{A}_{2,1}[2]$ and 
the product space $\mathcal{M}_2[2] \times \A_1[2]$.
Let $C$ be a genus $2$ curve with level $2$ structure represented
by the symplectic basis $(e_1,e_2,f_1,f_2)$ of $\mathrm{Jac}(C)[2]$ 
and let $E$ be an elliptic curve with level
$2$ structure $(e_3,f_3)$. Then $t_2(C) \times E$ is an abelian threefold and
$t_2(C)[2] \times E[2]$ is a six dimensional vector space over $\mathbb{F}_2$ with a symplectic pairing given by
\begin{equation*}
 e_i \cdot e_j = f_i \cdot f_j = 0
\end{equation*}
and
\begin{equation*}
 e_i \cdot f_j = \delta_{i,j}
\end{equation*}
for all $i$ and $j$, where we identify $e_i$ with $t_2(e_i)$ and $f_i$ with $t_2(f_i)$.
Clearly, not all level $2$ structures on $t_2(C) \times E$ arise in this way
but those that do are permuted by the group $\Sp{4} \times \Sp{2}$. Let $\mathscr{C}$ be the
quotient set $\Sp{6}/(\Sp{4} \times \Sp{2})$. We may then describe the locus $\mathcal{A}_{2,1}[2]$ as
\begin{equation*}
 \mathcal{A}_{2,1}[2] \cong \coprod_{c \in \mathscr{C}} (\mathcal{M}_2[2] \times \mathcal{A}_1[2])_c,
\end{equation*}
where $(\mathcal{M}_2[2] \times \mathcal{A}_1[2])_c$ is an isomorphic copy of $\mathcal{M}_2[2] \times \mathcal{A}_1[2]$ indexed by $c$ and the
components are permuted as
\begin{equation*}
 g(\mathcal{M}_2[2] \times \mathcal{A}_1[2])_c = (\mathcal{M}_2[2] \times \mathcal{A}_1[2])_{gc}
\end{equation*}
for $g \in \Sp{6}$. In terms of cohomology groups this means that
\begin{equation*}
 H^i(\mathcal{A}_{2,1}[2]) = \mathrm{Ind}_{\Sp{4} \times \Sp{2}}^{\Sp{6}} H^i(\mathcal{M}_2[2] \times \mathcal{A}_1[2]).
\end{equation*}
By the K\"unneth theorem we have that
\begin{equation*}
 H^i(\mathcal{M}_2[2] \times \mathcal{A}_1[2]) \cong \bigoplus_{p+q=i} H^p(\mathcal{M}_2[2]) \otimes H^q(\mathcal{A}_1[2]).
\end{equation*}
Thus, in order to understand the action of $\Sp{4} \times \Sp{2}$ on $H^i(\mathcal{M}_2[2] \times \mathcal{A}_1[2])$ it is
enough to understand the action of $\Sp{4}$ on $H^i(\mathcal{M}_2[2])$ and the action of $\Sp{2}$ on $H^i(\mathcal{A}_1[2])$ for all $i$.

\subsection{The moduli space of elliptic curves with level two structure} \label{sec-g1}
In order to understand the action of $\Sp{2}$ on $H^i(\mathcal{A}_1[2])$ we note that
$\Sp{2}$ is isomorphic to the symmetric group $S_3$ and that $\mathcal{A}_1[2]$ is isomorphic
to $\mathcal{M}_{0,4}$, the moduli space of four ordered points on $\mathbb{P}^1$. Under these
identifications, the action of $\Sp{2}$ is given by permuting the first three points. 

Since
$\mathcal{M}_{0,4}$ is pure in the sense of Dimca and Lehrer \cite{dimcalehrer} we can 
deduce the action of $\Sp{2}$ on the cohomology groups by 
a twisted point count, see Section~\ref{sec-euler}. 
Simple computations give, where $F$ is the Frobenius, 
\begin{equation*}
\resizebox{0.55\hsize}{!}{$
\arraycolsep=8pt\def\arraystretch{1.5}
\begin{array}{lll}
 |\mathcal{M}_{0,4}^{F \circ \mathrm{id}}|  & = \frac{(q+1)  q  (q-1)  (q-2)}{|\mathrm{PGL}_2(\mathbb{F}_q)|} & = q-2 \\
 |\mathcal{M}_{0,4}^{F \circ (12)}| & = \frac{(q+1)  (q^2-q)  q}{|\mathrm{PGL}_2(\mathbb{F}_q)|} & = q \\
 |\mathcal{M}_{0,4}^{F \circ (123)}| & = \frac{(q+1)  (q^3-q)}{|\mathrm{PGL}_2(\mathbb{F}_q)|}   & = q+1.
\end{array}$}
\end{equation*}
Thus, the traces of $(\mathrm{id},(12),(123))$ on 
$H^0(\mathcal{A}_1[2])$ and $H^1(\mathcal{A}_1[2])$ are
$(1,1,1)$ and $(2,0,-1)$, respectively. In other words, $H^0(\mathcal{A}_1[2])$ is the trivial representation
of $\Sp{2}$ while $H^1(\mathcal{A}_1[2])$ is the standard representation.

\subsection{The moduli space of genus two curves with level two structure} \label{sec-g2}
In order to understand the action of $\Sp{4}$ on $H^i(\mathcal{M}_2[2])$ we note that
$\Sp{4}$ is isomorphic to the symmetric group $S_6$ and that $\mathcal{M}_2[2]$ is isomorphic
to $\mathcal{M}_{0,6}$, the moduli space of six ordered points on $\mathbb{P}^1$.
Under these identifications,
the action of $\Sp{4}$ on $\mathcal{M}_2[2]$ is given by permuting the points.
Also, $\mathcal{M}_{0,6}$ is pure so we can again deduce the action of $\Sp{4}$ on the cohomology
via twisted point counts. Simple computations give, where $F$ is the Frobenius, 
\begin{equation*}
\resizebox{0.83\hsize}{!}{$
\arraycolsep=8pt\def\arraystretch{1.5}
\begin{array}{lll}
 |\mathcal{M}_{0,6}^{F \circ \mathrm{id}}| & = \frac{(q+1)  q  (q-1)  (q-2)  (q-3)  (q-4)}{|\mathrm{PGL}_2(\mathbb{F}_q)|}  & = q^3 - 9q^2 + 26q - 24\\ 
 |\mathcal{M}_{0,6}^{F \circ (12)}| & = \frac{(q+1)  q  (q-1)  (q-2)  (q^2-q)}{|\mathrm{PGL}_2(\mathbb{F}_q)|}  & = q^3 - 3q^2 + 2q \\
 |\mathcal{M}_{0,6}^{F \circ (12)(34)}| & = \frac{(q+1)  q  (q^2-q)  (q^2-q-2)}{|\mathrm{PGL}_2(\mathbb{F}_q)|}  & = q^3 - q^2 - 2q\\
 |\mathcal{M}_{0,6}^{F \circ(12)(34)(56)}| &  = \frac{(q^2-q)  (q^2-q-2)  (q^2-q-4)}{|\mathrm{PGL}_2(\mathbb{F}_q)|}  & = q^3 - 3q^2 - 2q + 8 \\
 |\mathcal{M}_{0,6}^{F \circ (123)}| & = \frac{(q+1)  q  (q-1)  (q^3-q)}{|\mathrm{PGL}_2(\mathbb{F}_q)|} & = q^3-q \\
 |\mathcal{M}_{0,6}^{F \circ (123)(45)}| & = \frac{(q+1) (q^2-q)  (q^3-q)}{|\mathrm{PGL}_2(\mathbb{F}_q)|} & = q^3-q \\
 |\mathcal{M}_{0,6}^{F \circ (123)(456)}| & = \frac{(q^3-q)  (q^3-q-3)}{|\mathrm{PGL}_2(\mathbb{F}_q)|} & = q^3-q-3 \\
 |\mathcal{M}_{0,6}^{F \circ (1234)}| & = \frac{(q+1)  q  (q^4-q^2)}{|\mathrm{PGL}_2(\mathbb{F}_q)|} & = q^3+q^2 \\
 |\mathcal{M}_{0,6}^{F \circ (1234)(56)}| & = \frac{(q^2-q)  q  (q^4-q^2)}{|\mathrm{PGL}_2(\mathbb{F}_q)|} & = q^3 - q^2 \\
 |\mathcal{M}_{0,6}^{F \circ (12345)}| & = \frac{(q+1)  (q^5-q)}{|\mathrm{PGL}_2(\mathbb{F}_q)|} & = q^3+q^2+q+1 \\
|\mathcal{M}_{0,6}^{F \circ (123456)}| & = \frac{q^6-q^3-q^2+q}{|\mathrm{PGL}_2(\mathbb{F}_q)|} & = q^3+q-1.
\end{array}$}
\end{equation*}
We may read off $\mathrm{Tr}(\sigma,H^i(\mathcal{M}_{0,6}))$ as the coefficient of $(-1)^i  q^{3-i}$ in
$|\mathcal{M}_{0,6}^{F \circ \sigma'}|$, where $\sigma'$ is the element in the above list which is conjugate to $\sigma$.

We now know $H^i(\mathcal{M}_2[2])$ as a representation of $\Sp{4}$ for all $i$ and we know
$H^i(\mathcal{A}_1[2])$ as a representation of $\Sp{2}$ for all $i$ so we therefore know
$H^{i}(\mathcal{M}_2[2] \times \mathcal{A}_1[2])$ as a representation of $\Sp{4} \times \Sp{2}$.
Inducing this representation to $\Sp{6}$ gives us $H^i(\mathcal{A}_{2,1}[2])$ as a representation of
$\Sp{6}$. We give the result in
Table~\ref{a21table}. 

As an aside we note that the complement of $t_2(\mathcal{M}_{2}[2])$ 
inside $\mathcal{A}_{2}[2]$ consists of products of elliptic curves. The cohomology  of this 
locus can be computed using the same techniques as in Section~\ref{sec-A111}, but in this 
case we will
omit the details. Let us denote  
the irreducible representations of $\Sp{4} \cong  S_6$ by $s_{\lambda}$, which are indexed 
in the standard way by $\lambda$, a partition of~$6$. Adding the contributions from the two loci we find that,  
$$ e_{\mathcal{A}_2[2]}(v) =(1+v^2)s_{6} -v^4(s_{5,1}+s_{4,2})+v^6s_{3,2,1}.
$$

\section{Products of three elliptic curves} \label{sec-A111} 
There is a close relationship between the product space $(\A_1[2])^3$ and the locus
$\mathcal{A}_{1,1,1}[2]$. Let $E_1$, $E_2$ and $E_3$ be three elliptic curves
with level $2$ structures $(e_{1},f_{1})$, $(e_{2},f_2)$ and $(e_3,f_3)$, respectively.
Then $E_1 \times E_2 \times E_3$ is an abelian threefold and $E_1[2] \times E_2[2] \times E_3[2]$
is a six dimensional vector space over $\mathbb{F}_2$ with a symplectic pairing given by
\begin{equation*}
 e_i \cdot e_j = f_i \cdot f_j = 0
\end{equation*}
for all $i$ and $j$ and
\begin{equation*}
 e_i \cdot f_j = \delta_{i,j}.
\end{equation*}
Clearly, not all level $2$ structures on $E_1 \times E_2 \times E_3$ arise in this way
but those that do are permuted by the group $(\Sp{2})^3$ while the three curves themselves
are permuted by the symmetric group $S_3$. Let $\mathscr{C}$ be the
quotient set $\Sp{6}/(S_3 \ltimes (\Sp{2})^3)$. We may describe the locus $\mathcal{A}_{1,1,1}[2]$ as
\begin{equation*}
 \mathcal{A}_{1,1,1}[2] \cong \coprod_{c \in \mathscr{C}} (\mathcal{A}_1[2])^3_c,
\end{equation*}
where $(\mathcal{A}_1[2])^3_c$ is an isomorphic copy of $(\mathcal{A}_1[2])^3$ indexed by $c$ and the
components are permuted as
\begin{equation*}
 g(\mathcal{A}_1[2])^3_c = (\mathcal{A}_1[2])^3_{gc}
\end{equation*}
for $g \in \Sp{6}$. In terms of cohomology groups this means that
\begin{equation*}
 H^i(\mathcal{A}_{1,1,1}[2]) = \mathrm{Ind}_{S_3 \ltimes (\Sp{2})^3}^{\Sp{6}} H^i((\mathcal{A}_1[2])^3).
\end{equation*}
By the K\"unneth theorem we have that
\begin{equation*}
 H^i((\mathcal{A}_1[2])^3) \cong \bigoplus_{p+q+r=i} H^p(\mathcal{A}_1[2]) \otimes H^q(\mathcal{A}_1[2]) \otimes H^r(\mathcal{A}_1[2]).
\end{equation*}
Thus, in order to understand the action of $S_3 \ltimes (\Sp{2})^3$ on $H^i((\mathcal{A}_1[2])^3$ it is
enough to understand the action of $\Sp{2}$ on $H^i(\mathcal{A}_1[2])$ for all $i$ and the action of $S_3$ on the 
factors.

Since the action of $\Sp{2}$ was described in Section~\ref{sec-g1} we only consider the action of $S_3$ on the factors.
Let $\alpha \otimes \beta \otimes \gamma \in H^p(\mathcal{A}_1[2]) \otimes H^q(\mathcal{A}_1[2]) \otimes H^r(\mathcal{A}_1[2]) \subseteq H^{p+q+r}((\mathcal{A}_1[2])^3)$.
We have
\begin{align*}
 (12).\alpha \otimes \beta \otimes \gamma = (-1)^{pq}\beta \otimes \alpha \otimes \gamma, \\
 (13).\alpha \otimes \beta \otimes \gamma = (-1)^{pr}\gamma \otimes \beta \otimes \alpha, \\
 (23).\alpha \otimes \beta \otimes \gamma = (-1)^{qr}\alpha \otimes \gamma \otimes \beta,
\end{align*}
where the signs are a consequence of the K\"unneth isomorphism.
Since $S_3$ is generated by transpositions and $H^{p+q+r}((\mathcal{A}_1[2])^3)$ is generated by
elements of the form $\alpha \otimes \beta \otimes \gamma$ for all possible choices of $p$, $q$ and
$r$, this determines the action of $S_3$ on $H^{p+q+r}((\mathcal{A}_1[2])^3)$.

We now have all the information we need in order to understand the action of $S_3 \ltimes (\Sp{2})^3$ on $H^i((\mathcal{A}_1[2])^3)$.

\pagebreak[2]
\begin{example}
\label{a1example}
 Let $u$ be a basis vector for the trivial representation of $\Sp{2}$ and
 let $v_1$ and $v_2$ be basis vectors for the standard representation of $\Sp{2}$.
 Let $\sigma \in \Sp{2}$ be an element of order $3$ acting as
 \begin{align*}
  \sigma.u & = u, \\
  \sigma.v_1 & = v_2, \\
  \sigma.v_2 & = -v_1-v_2,
 \end{align*}
and let $g \in S_3 \ltimes (\Sp{2})^3$ on $H^i((\mathcal{A}_1[2])^3$ be the element $g=((23),(\sigma,\sigma,\mathrm{id}))$.
We of course have $g.u\otimes u \otimes u = u \otimes u \otimes u$, so $\mathrm{Tr}(g,H^0((\mathcal{A}_1[2])^3)=1$.
In cohomological degree $1$ we have
\begin{equation*}
 g.v_2 \otimes u \otimes u = -(v_1 + v_2) \otimes u \otimes u
\end{equation*}
while $g.\alpha$ has no component in the direction of $\alpha$ for all other choices of $\alpha \in H^1((\mathcal{A}_1[2])^3)$.
Thus, $\mathrm{Tr}(g,H^1((\mathcal{A}_1[2])^3)=-1$. In degree $2$ we have
\begin{equation*}
 g.u \otimes v_2 \otimes v_2 = u \otimes v_2 \otimes (v_1+v_2)
\end{equation*}
while $g.\alpha$ has no component in the direction of $\alpha$ for all other choices of $\alpha \in H^2((\mathcal{A}_1[2])^3)$.
We conclude that $\mathrm{Tr}(g,H^2((\mathcal{A}_1[2])^3)=1$. Finally, in degree $3$ we have
\begin{equation*}
 g.v_2 \otimes v_2 \otimes v_2 = -(v_1+v_2) \otimes v_2 \otimes (v_1+v_2)
\end{equation*}
while $g.\alpha$ has no component in the direction of $\alpha$ for all other choices of $\alpha \in H^3((\mathcal{A}_1[2])^3)$.
Hence, $\mathrm{Tr}(g,H^3((\mathcal{A}_1[2])^3)=-1$. We thus have
\begin{equation*}
 \sum_{i=0}^3 \mathrm{Tr}(g,H^i((\mathcal{A}_1[2])^3) \, t^i = 1-t+t^2-t^3. 
\end{equation*}
\end{example}

Similar computations for the other conjugacy classes of $S_3 \ltimes (\Sp{2})^3$ give
the results in Table~\ref{a1table}, where
$$P_g((\mathcal{A}_1[2])^3,t):=\sum_{i=0}^3 \mathrm{Tr}(g,H^i((\mathcal{A}_1[2])^3)) \, t^i,$$
is called the equivariant Poincar\'e polynomial.
See Chapter 4 of \cite{jameskerber} for a beautiful description
of how to compute representatives of $S_3 \ltimes (\Sp{2})^3$. In Table~\ref{a1table}, $\sigma$ is
the element of $\Sp{2}$ described in Example~\ref{a1example} while $\tau$ is the element of order $2$
acting as
\begin{equation*}
 \tau.v_1=-v_1, \quad \tau.v_2=v_1+v_2
\end{equation*}
where $v_1$ and $v_2$ are the same basis vectors of the standard representation considered in Example~\ref{a1example}.

\begin{table}[htbp]
\begin{equation*}
\resizebox{0.95\textwidth}{!}{$
 \begin{array}{l|l|l|l|l|l}
  g & P_g((\mathcal{A}_1[2])^3,t) & g & P_g((\mathcal{A}_1[2])^3,t) & g & P_g((\mathcal{A}_1[2])^3,t) \\
  \hline
  (\textrm{id},(\textrm{id},\textrm{id},\textrm{id})) & 1 + 6t+12t^2+8t^3 & ((23),(\textrm{id},\textrm{id},\textrm{id})) & 1+2t-2t^2-4t^3 & ((123),(\textrm{id},\textrm{id},\textrm{id})) & 1+2t^3 \\
  (\textrm{id},(\textrm{id},\textrm{id},\tau)) & 1+4t+4t^2 & ((23),(\textrm{id},\tau,\textrm{id})) & 1+2t &
  ((123),(\tau,\textrm{id},\textrm{id})) & 1 \\
  (\textrm{id},(\textrm{id},\textrm{id},\sigma)) & 1+3t-4t^3 & ((23),(\textrm{id},\sigma,\textrm{id})) & 1+2t+t^2+2t^3 & ((123),(\sigma,\textrm{id},\textrm{id})) & 1-t^3 \\
  (\textrm{id},(\textrm{id},\tau,\tau)) & 1+2t & ((23),(\tau,\textrm{id},\textrm{id})) & 1-2t^2 & \, & \, \\
  (\textrm{id},(\textrm{id},\tau,\sigma)) & 1+t-2t^2 & ((23),(\tau,\tau,\textrm{id})) & 1 & \, & \,\\
  (\textrm{id},(\textrm{id},\sigma,\sigma)) & 1 -3t^2+2t^3 & ((23),(\tau,\sigma,\textrm{id})) & 1+t^2 & \, & \, \\
  (\textrm{id},(\tau,\tau,\tau)) & 1 & ((23),(\sigma,\textrm{id},\textrm{id})) & 1-t-2t^2+2t^3 & \, & \,\\
  (\textrm{id},(\tau,\tau,\sigma)) & 1-t & ((23),(\sigma,\tau,\textrm{id})) & 1-t & \, & \,\\
  (\textrm{id},(\tau,\sigma,\sigma)) & 1-2t+t^2 & ((23),(\sigma,\sigma,\textrm{id})) & 1-t+t^2-t^3 & \, & \,\\
  (\textrm{id},(\sigma,\sigma,\sigma)) & 1-3t+3t^2-t^3 & \, & \,& \, & \,
 \end{array}$}
\end{equation*}
\caption{Equivariant Poincar\'e polynomials of $(\mathcal{A}_1[2])^3$ for a representative $g$ of every conjugacy class of $S_3 \ltimes (\Sp{2})^3$.}
\label{a1table}
\end{table}

By inducing the corresponding representations from  $S_3 \ltimes (\Sp{2})^3$  to $\Sp{6}$ we obtain 
the cohomology of $\mathcal{A}_{1,1,1}[2]$ as a representation of $\Sp{6}$. We give the result in
Table~\ref{a111table}. 

\section{Cohomology groups of strata}

In this section we give the cohomology groups of $\mathcal{Q}[2]$, $\mathcal{H}_3[2]$,
$\mathcal{A}_{2,1}[2]$ and $\mathcal{A}_{1,1,1}[2]$ as representations of $\Sp{6}$.
The results are presented in Table~\ref{Qtable}-\ref{a111table}.
Each column in these tables corresponds to an irreducible representation of $\Sp{6}$.
The irreducible representations are denoted $\phi_{dn}$ where $d$ is the dimension of
the representation and $n$ is letter used to distinguish different representations of the same
dimension, see~\cite{conwayetal}.

\begin{table}[htbp]
\begin{equation*}
\resizebox{0.62\textwidth}{!}{$
\begin{array}{r|rrrrrrrrrr} 
\, & \phi_{1a} & \phi_{7a} & \phi_{15a} & \phi_{21a} & \phi_{21b} & \phi_{27a} & \phi_{35a} & \phi_{35b} & \phi_{56a} & \phi_{70a} \\
\hline
H^0 & 1 & 0 & 0 & 0 & 0 & 0 & 0 & 0 & 0 & 0 \\
H^1 & 0 & 0 & 0 & 0 & 0 & 0 & 0 & 1 & 0 & 0 \\
H^2 & 0 & 0 & 0 & 0 & 0 & 0 & 0 & 0 & 0 & 0 \\
H^3 & 0 & 0 & 0 & 1 & 0 & 0 & 0 & 0 & 0 & 0 \\
H^4 & 0 & 0 & 0 & 0 & 0 & 0 & 0 & 0 & 0 & 1 \\
H^5 & 0 & 0 & 0 & 0 & 0 & 1 & 1 & 1 & 0 & 0 \\
H^6 & 1 & 0 & 2 & 0 & 1 & 1 & 1 & 3 & 0 & 0 \\
\hline
\, & \phi_{84a} & \phi_{105a} & \phi_{105b} & \phi_{105c} & \phi_{120a} & \phi_{168a} & \phi_{189a} & \phi_{189b} & \phi_{189c} & \phi_{210a} \\
\hline
H^0 & 0 & 0 & 0 & 0 & 0 & 0 & 0 & 0 & 0 & 0 \\
H^1 & 0 & 0 & 0 & 0 & 0 & 0 & 0 & 0 & 0 & 0 \\
H^2 & 0 & 0 & 0 & 0 & 0 & 0 & 0 & 0 & 0 & 1 \\
H^3 & 0 & 0 & 1 & 0 & 0 & 0 & 1 & 0 & 0 & 2 \\
H^4 & 0 & 0 & 2 & 0 & 2 & 1 & 2 & 1 & 0 & 3 \\
H^5 & 1 & 2 & 2 & 1 & 2 & 4 & 3 & 3 & 3 & 4 \\
H^6 & 5 & 1 & 1 & 4 & 0 & 3 & 2 & 2 & 5 & 3 \\
\hline
\, & \phi_{210b} & \phi_{216a} & \phi_{280a} & \phi_{280b} & \phi_{315a} & \phi_{336a} & \phi_{378a} & \phi_{405a} & \phi_{420a} & \phi_{512a} \\
\hline
H^0 & 0 & 0 & 0 & 0 & 0 & 0 & 0 & 0 & 0 & 0 \\
H^1 & 0 & 0 & 0 & 0 & 0 & 0 & 0 & 0 & 0 & 0 \\
H^2 & 0 & 0 & 0 & 1 & 0 & 0 & 0 & 0 & 0 & 0 \\
H^3 & 1 & 0 & 0 & 0 & 0 & 0 & 1 & 2 & 2 & 1 \\
H^4 & 4 & 0 & 3 & 1 & 3 & 2 & 3 & 6 & 5 & 4 \\
H^5 & 4 & 4 & 4 & 6 & 5 & 6 & 6 & 6 & 8 & 9 \\
H^6 & 1 & 6 & 3 & 6 & 1 & 6 & 4 & 2 & 6 & 6
\end{array}$}
\end{equation*}
\caption{The cohomology groups of $\mathcal{Q}[2]$ as a representation of $\Sp{6}$.}
\label{Qtable}
\end{table}

 \begin{table}[htbp]
\begin{equation*}
\resizebox{0.62\textwidth}{!}{$
\begin{array}{r|rrrrrrrrrr} 
\, & \phi_{1a} & \phi_{7a} & \phi_{15a} & \phi_{21a} & \phi_{21b} & \phi_{27a} & \phi_{35a} & \phi_{35b} & \phi_{56a} & \phi_{70a} \\
\hline
H^0 & 1&0&0&0&0&0&0&1&0&0\\ 
H^1 & 0&0&0&0&0&1&0&1&0&0\\ 
H^2 & 0&0&0&1&0&0&0&0&0&0\\ 
H^3 & 0&0&0&1&0&0&0&0&0&1\\ 
H^4 & 0&0&0&0&0&1&1&1&0&1\\ 
H^5 & 0&0&1&0&1&1&1&2&0&0 \\
\hline
\, & \phi_{84a} & \phi_{105a} & \phi_{105b} & \phi_{105c} & \phi_{120a} & \phi_{168a} & \phi_{189a} & \phi_{189b} & \phi_{189c} & \phi_{210a} \\
\hline
H^0 & 0&0&0&0&0&0&0&0&0&0\\ 
H^1 & 0&0&0&0&0&1&0&0&0&1\\
H^2 & 0&0&1&0&2&1&1&0&0&3\\
H^3 & 0&0&3&1&3&2&4&1&0&5\\
H^4 & 2&2&3&2&3&6&5&4&4&6\\
H^5 & 4&2&1&4&1&4&3&3&6&4 \\
\hline
\, & \phi_{210b} & \phi_{216a} & \phi_{280a} & \phi_{280b} & \phi_{315a} & \phi_{336a} & \phi_{378a} & \phi_{405a} & \phi_{420a} & \phi_{512a} \\
\hline
H^0 & 0&0&0&0&0&0&0&0&0&0\\
H^1 & 0&0&0&1&0&0&0&0&0&0\\
H^2 & 1&0&0&2&0&0&1&3&2&2\\
H^3 & 5&1&3&3&4&3&5&10&7&7\\
H^4 & 6&5&7&8&7&9&9&10&12&14\\
H^5 & 2&7&4&8&3&8&6&4&8&9
\end{array}$}
\end{equation*}
\caption{The cohomology groups of $\mathcal{H}_{3}[2]$ as a representation of $\Sp{6}$.}
\label{Hyptable}
\end{table}

\begin{table}[ht]
 \begin{equation*}
 \resizebox{0.62\textwidth}{!}{$
  \begin{array}{c|cccccccccc}
  \, & \phi_{1a} & \phi_{7a} & \phi_{15a} & \phi_{21a} & \phi_{21b} & \phi_{27a} & \phi_{35a} & \phi_{35b} & \phi_{56a} & \phi_{70a} \\
  \hline
  H^0 & 1&0&0&0&0&1&0&1&0&0 \\
  H^1 & 0&0&0&0&0&2&0&2&0&0 \\
  H^2 & 0&0&0&1&0&0&0&1&0&0 \\
  H^3 & 0&0&0&1&0&0&0&0&0&2 \\
  H^4 & 0&0&0&0&0&0&0&0&1&1 \\
  \hline
  \, & \phi_{84a} & \phi_{105a} & \phi_{105b} & \phi_{105c} & \phi_{120a} & \phi_{168a} & \phi_{189a} & \phi_{189b} & \phi_{189c} & \phi_{210a} \\
  \hline
  H^0 & 0&0&1&0&0&1&0&0&0&0 \\
  H^1 & 1&0&1&0&2&2&0&0&0&2 \\
  H^2 & 1&0&2&2&3&3&3&0&0&5 \\ 
  H^3 & 1&0&3&2&2&3&5&2&1&6 \\
  H^4 & 0&1&1&0&2&1&2&2&1&2 \\
  \hline
  \, & \phi_{210b} & \phi_{216a} & \phi_{280a} & \phi_{280b} & \phi_{315a} & \phi_{336a} & \phi_{378a} & \phi_{405a} & \phi_{420a} & \phi_{512a} \\
  \hline
  H^0 & 0&0&0&0&0&0&0&0&0&0 \\
  H^1 & 1&0&0&3&0&0&0&1&1&1 \\
  H^2 & 2&2&0&6&1&2&3&6&5&5 \\
  H^3 & 4&2&4&4&4&5&6&10&8&10 \\
  H^4 & 4&1&3&2&5&3&5&6&5&5
  \end{array}$}
 \end{equation*}
 \caption{The cohomology groups of $\mathcal{A}_{2,1}[2]$ as representations of $\Sp{6}$.}
 \label{a21table}
\end{table}

\begin{table}[htbp]
\begin{equation*}
\resizebox{0.62\textwidth}{!}{$
 \begin{array}{c|cccccccccc}
 \, & \phi_{1a} & \phi_{7a} & \phi_{15a} & \phi_{21a} & \phi_{21b} & \phi_{27a} & \phi_{35a} & \phi_{35b} & \phi_{56a} & \phi_{70a} \\
 \hline
 H^0 & 1&0&0&0&0&1&0&1&0&0 \\
 H^1 & 0&0&0&0&0&1&0&2&0&0 \\
 H^2 & 0&0&0&1&0&0&0&0&0&1 \\
 H^3 & 0&0&0&0&0&0&0&0&1&1 \\
 \hline
 \, & \phi_{84a} & \phi_{105a} & \phi_{105b} & \phi_{105c} & \phi_{120a} & \phi_{168a} & \phi_{189a} & \phi_{189b} & \phi_{189c} & \phi_{210a} \\
 \hline
 H^0 & 1&0&1&0&0&1&0&0&0&0 \\
 H^1 & 1&0&1&1&2&2&1&0&0&3 \\
 H^2 & 0&0&2&1&2&2&3&1&0&4 \\
 H^3 & 0&0&1&0&1&0&1&1&0&1 \\
 \hline
 \, & \phi_{210b} & \phi_{216a} & \phi_{280a} & \phi_{280b} & \phi_{315a} & \phi_{336a} & \phi_{378a} & \phi_{405a} & \phi_{420a} & \phi_{512a} \\
 \hline
 H^0 & 0&0&0&1&0&0&0&0&1&0 \\
 H^1 & 1&1&0&4&0&1&1&2&2&2 \\
 H^2 & 2&1&1&3&2&2&3&6&4&5 \\
 H^3 & 3&0&2&0&3&1&3&4&3&3
 \end{array}$}
\end{equation*}
\caption{The cohomology groups of $\mathcal{A}_{1,1,1}[2]$ as representations of $\Sp{6}$.}
\label{a111table}
\end{table}

  \clearpage
\bibliographystyle{plain}

\renewcommand{\bibname}{References}

\bibliography{references} 
\end{document}